\documentclass[fleqn]{article}
\usepackage[english]{babel}
\usepackage{a4wide} 
\usepackage{mathrsfs}
\usepackage{amsmath}
\usepackage{amsfonts}
\usepackage{amssymb}
\usepackage{amsthm}
\usepackage{verbatim}
\usepackage{graphicx}
\usepackage[colorlinks=false,pdfborderstyle={/W 1}]{hyperref}
\def\arXiv#1#2{\href{http://front.math.ucdavis.edu/#1}{{\tt arXiv:#1 [#2]}}}

\input epsf.sty

\newtheorem{thm}{Theorem}[section]
\newtheorem{quest}[thm]{Question}

\newtheorem{lemma}[thm]{Lemma}
\newtheorem{proposition}[thm]{Proposition}
\newtheorem{corollary}[thm]{Corollary}
\theoremstyle{remark}\newtheorem{rem}[thm]{Remark}
\theoremstyle{definition}
\newtheorem{definition}{Definition}[section]

\newcommand{\R}{\mathbb{R}}

\newcommand{\E}{\mathbb{E}}

\newcommand{\N}{\mathbb{N}}

\newcommand{\Maj}{\mathsf{Maj}}
\newcommand{\Piv}{\mathsf{Piv}}

\numberwithin{equation}{section}
\numberwithin{figure}{section}

\def\bl{\begin{lemma}}
\def\el{\end{lemma}}
\def\bth{\begin{thm}}
\def\eth{\end{thm}}
\def\bc{\begin{corollary}}
\def\ec{\end{corollary}}
\def\bcj{\begin{conjecture}}
\def\ecj{\end{conjecture}}
\def\bpr{\begin{proposition}}
\def\epr{\end{proposition}}
\def\bde{\begin{definition}}
\def\ede{\end{definition}}

\def\beq{\begin{equation}}
\def\eeq{\end{equation}}

\def\bpf{\begin{proof}}
\def\epf{\end{proof}}

\def \eps {\epsilon}

\def\bl{\begin{lemma}}
\def\el{\end{lemma}}
\def\bc{\begin{corollary}}
\def\ec{\end{corollary}}
\def\bcj{\begin{conjecture}}
\def\ecj{\end{conjecture}}
\def\bpr{\begin{proposition}}
\def\epr{\end{proposition}}
\def\bde{\begin{definition}}
\def\ede{\end{definition}}
\newcommand{\be}{\begin{eqnarray}}
\newcommand{\ee}{\end{eqnarray}}
\newcommand{\bes}{\begin{eqnarray*}}
\newcommand{\ees}{\end{eqnarray*}}
\newcommand{\comm}[1]{}
\usepackage{mathrsfs}

\def\E{\mathbb{E}}

\def\N{\mathbb{N}}
\def\R{\mathbb{R}}
\def\Pr{\mathbb{P}}

\def\Piv{\mathscr{P}}

\def\eps{\epsilon}
\def \Maj{\mathsf{Maj}}
\def \Tr{\mathsf{Tribes}}

\def\1{1\!\! 1}

\begin{document}

\title{Pivotality versus noise
stability for monotone transitive functions}
\author{P\'al Galicza \\
Central European University\\
Alfr\'ed R\'enyi Institute of Mathematics}

\date{}
\maketitle
\begin{abstract}
 We construct a noise stable sequence of transitive, monotone increasing Boolean functions  $f_n: \{-1,1\}^{k_n} \longrightarrow \{-1,1\}$ which admit many pivotals with high probability. We show that such a sequence is volatile as well, and thus it is also an example of a volatile and noise stable sequence of transitive, monotone functions. 
\end{abstract}

\section{Introduction}

In this note we answer the following question posed by Gil Kalai and Gady Kozma (Oberwolfach, September 2018):
 Is there a sequence of Boolean functions $f_n: \left\{ -1,1 \right\}^{k_n} \longrightarrow  \left\{-1, 1 \right\}$ such that $f_n$ is transitive and noise stable, but at the same time $\Pr[\Piv_n(\omega)\neq \emptyset] > c$ for some constant $c>0$ for all $n \in \N$?

We are going to show that the answer is positive.
\bth \label{pivstable}
There exists a sequence of transitive monotone functions $f_n: \left\{ -1,1 \right\}^{k_n} \longrightarrow  \left\{-1, 1 \right\}$ such that $f_n$ is transitive, noise stable and $\lim_{n}{\Pr[\Piv_n > a_n]} =1 $ (here $\Piv_n$ is the pivotal set, see below) for some sequence of integers $a_n \rightarrow \infty $.
\eth

We begin by giving the necessary definitions. For further details on noise sensitivity, noise stability and their relation to the pivotal set see \cite{GS}.

We call a function $f: \left\{ -1,1 \right\}^{n} \longrightarrow  \R$ transitive if there is a transitive group action on the coordinate set $k$ and $f$ is invariant under this action.

\bde [Pivotal Set]
The pivotal set for a Boolean function $f: \left\{ -1,1 \right\}^{n} \longrightarrow  \left\{-1, 1 \right\}$ is the random set of coordinates $i$ for which  $f(\omega) \neq f(\omega^i)$, where $\omega^i$ is $\omega$ with its $i$th coordinate flipped.
\ede

\begin{definition} [Noise Sensitivity]
Let $\eps$ be a  positive real number. For a uniform random vector $\omega \in \{-1,1\}^{k_n} $ denote $N_{\epsilon}(\omega)$ the random vector which we obtain from $\omega$ by resampling each of its bits independently with probability $\epsilon$. A sequence of functions $f_n: \{-1,1\}^{k_n} \longrightarrow \{-1,1\}$  is noise sensitive if 
\begin{equation} 
\displaystyle \lim_{n\rightarrow \infty} \text{Corr}(f_n (\omega), f_n (N_{\eps}(\omega)))= 0  
\end{equation}
\end {definition}

\begin{definition} [Noise Stability]
 A sequence of functions $f_n: \{-1,1\}^{k_n} \longrightarrow \{-1,1\}$  is noise stable if 
\begin{equation} 
\displaystyle \lim_{\eps \rightarrow 0} \sup_{n} \Pr[f_n (\omega) \neq f_n (N_{\eps}(\omega))]= 0  
\end{equation}
\end {definition}

In \cite{GS} Section XII. 2 a sequence of monotone, non-degenerate and noise stable Boolean functions is presented which has many pivotals with a positive probability (The construction is due to O. Schramm). This example, however, is not transitive.

These types of questions are interesting since noise sensitivity and having many pivotals are often closely related to each other; see, e.g., the case of crossing events in planar percolation \cite{GPS1}. Noise sensitivity and noise stability of a function can both be expressed in terms of the typical size of the Spectral Sample, which is also a random subset of the index set (see \cite{GS}). It is an interesting fact that the first and second order marginals of the pivotal set and the spectral sample are the same. This (falsely) suggests that for a noise stable functions the pivotal set is typically small. Our result is another indication that, in general, these two random sets may show very different behavior.  

Another dynamical property of Boolean functions, which may look, at first glance, almost the same as noise sensitivity, is volatility, studied in \cite{JS}. It roughly says that if we are updating the input bits in continuous time, then the output changes very often; see Definition \ref{volat} below. In Lemma \ref{piv_volatile} we show that having many pivotals implies being volatile, hence we obtain the following

\bc There exists a sequence of transitive monotone noise stable and volatile functions.  
\ec 
Our construction also implies, see Corollary \ref{dense_piv} below, that every (monotone) Boolean function is close to a (monotone) Boolean function that has many pivotals with high probability. As functions with these properties are also volatile, this is a strengthening of Theorem 1.4  in \cite{Fo}.

\comm{\begin{rem}
One can relax the stability condition to the lack of noise sensitivity.  In this case the answer is almost trivial.  We now sketch an example of a sequence of monotone functions which is transitive, not noise sensitive and the pivotal set is nonempty with a uniformly positive probability.

Let $A_n,B_n \subseteq \left\{ -1,1 \right\}^{k_n}$ be two sequences of monotone transitive events  satisfying the conditions of  Theorem \ref{pivstable} except for noise stability, with the property that there exists $c>0$ such that for all (large) $n$, $\Pr[ A_n \cup B_n] -\Pr[ A_n ]> c$ (We can, for example, choose  $A_n$ and $B_n$ to be two tribes events defined on different tribe partitions). Let $\Maj_{k_n}$ be the Majority function on the same $k_n$ bits. Now let 
$$
f_n =
\left\{
	\begin{array}{ll}
	\1_{A_n}	 & \mbox{if } \Maj_{k_n} = -1 \\
		\1_{A_n\cup B_n} & \mbox{if } \Maj_{k_n} = 1.
	\end{array}
\right.
$$   
It is clear that $f_n$ is monotone, transitive and admits pivotals with a positive probability. At the same time it is positively correlated with $\Maj_{k_n}$ and therefore cannot be noise sensitive.
\end{rem}}

\section{Construction}

In the sequel, we shall construct a sequence of functions $f_n: \left\{ -1,1 \right\}^{k_n} \longrightarrow  \left\{-1, 0, 1 \right\}$ with the following properties:

\begin{enumerate} 
\item $f_n$ is transitive
\item $ \displaystyle \lim _{n} \Pr[f_n  = 0] = 1$ 
\item $\displaystyle \lim _{n} \Pr[ \exists\; i, j \in  [k_n] \;: \;f_n(\omega^{i}) =1  \; \text{and}\;   f_n(\omega^{j}) = -1] = 1.$
\end{enumerate}
where $\omega^i$ denotes  $\omega$ with its $i$th coordinate flipped.
We will call a sequence of functions bribable if it satisfies the above  conditions.

Using a bribable sequence $f_n$ one can easily construct a transitive noise stable Boolean function which admits a pivotal bit with high probability. Namely, let $\Maj_n$ denote the majority function on the corresponding bit set. Let 
$$
g_n =
\left\{
	\begin{array}{ll}
		\Maj_n & \mbox{if } f_n = 0 \\
		f_n & \mbox{if } f_n \neq 0.
	\end{array}
\right.
$$
Obviously $g_n$ is noise stable because of property 2 of $f_n$. On the other hand, conditioned on $\left\{ f_n = 0 \right\} $ there is a pivotal bit with high probability because of property 3 of the sequence $f_n$.

It is also straightforward to verify that if we choose a bribable sequence $f_n$ which is monotone then the resulting $g_n$ sequence will be monotone as well.

Now we turn to the construction of a monotone bribable sequence. Define the Boolean function $\Tr(l,k): \left\{ -1,1 \right\}^{lk} \longrightarrow  \left\{0, 1 \right\}$ as follows: we group the bits in $k$ $l$-element subsets, these are the so called tribes. The function takes on $1$ if there is a tribe $T$ such that for every $i \in T \;:\;\omega(i)=1$, and $0$ otherwise. The $\Tr$ function is standard example, when $k_n$ and $l_n$ are defined in such a way that the function is non-degenerate. It is well know that such a sequence testifies that the Kahn-Kalai-Linial theorem about the maximal influence of sequences of Boolean functions (Theorem 1.14 in \cite{GS}) is sharp. 

We are going to show that in case the two sequences $l_n,k_n$ are properly chosen, a slight modification of $\Tr(l_n,k_n)$ is bribable.

\bpr
Suppose that $l_n$ and $k_n$ are sequences such that

\beq \label{kisVar}
\lim_{n \rightarrow \infty} {\left(1- \frac{1}{2^{l_n}} \right)^{k_n}} = 1
\eeq
and
\beq \label{nagyInf}
\lim_{n \rightarrow \infty} k_n l_n \frac{1}{2^{l_n}} = \infty
\eeq
then the sequence of functions $f_n (\omega):= \Tr(l_n,k_n)(\omega)- \Tr(l_n,k_n)(-\omega)$ is bribable. Moreover, there is a sequence of positive integers $a_n \rightarrow \infty$ such that $\Pr[ |\Piv_n| > a_n] \rightarrow 1$
\epr

\bpf

Let us call a tribe $T$ pivotal if there is exactly one $j \in T$ such that $\omega(j)=-1$. Define the random variable $X_n$ as the number of pivotal tribes in a configuration. Note that $\E[X_n] = k_n l_n\frac{1}{2^{l_n}}$. 

It is clear that conditioned on the event $\left\{ \Tr(l_n,k_n) =0 \right\}$ we have $|\Piv_n| =  X_n $, where $|\Piv_n|$ denotes the pivotal set of $\Tr(l_n,k_n)$. Consequently, for the respective conditional expected values:
$$
 \E [\Piv_n | \Tr(l_n,k_n) =0] = \E[X_n |\Tr(l_n,k_n) =0].
$$

We can write $X_n= \sum_{j}^{k_n}{Y_j}$ where $Y_j$ is the indicator of the event that the $j$th tribe is pivotal. For any $j \in [k_n]$ we have
$$
  \Pr[Y_j = 1 |\Tr(l_n,k_n) =1] =\frac{\Pr[Y_j = 1]  \Pr[\Tr(l_n,k_{n}-1) =1]}{\Pr[\Tr(l_n,k_n) =1]} \leq \Pr[Y_j = 1],
$$
using that if the  $j$th tribe is pivotal and there is a full $1$ tribe then the latter is among the remaining $k_{n}-1$  tribes.
This implies
$$
\E[X_n |\Tr(l_n,k_n) =1] \leq \E[X_n] \leq  \E[X_n |\Tr(l_n,k_n) =0] 
$$
and therefore 
$$
  \E[\Piv_n | \Tr(l_n,k_n) =0]  \geq   \E[X_n] = k_n l_n \frac{1}{2^{l_n}} \rightarrow \infty.
$$

As $X_n$ is binomially distributed with $\E[X_n]\rightarrow \infty $, being the sum of i.i.d $0-1$-valued random variables, there is a $a_n \rightarrow \infty$  such that 
$$
\lim_{n \rightarrow \infty} \Pr[X_n> a_n] =1.
$$
Note that
$$
\Pr[\Tr(l_n,k_n)= 0] = \left(1- \frac{1}{2^{l_n}} \right)^{k_n}
$$
and this probability tends to $1$ as $n$ approaches $\infty$ by our assumption. So clearly 
$$
\Pr[ X_n> a_n \; \text{and} \;  \Tr(l_n,k_n)= 0 ] = \Pr[ |\Piv_n| > a_n, \; \text{and} \;  \Tr(l_n,k_n) =0 ]  \rightarrow 1
$$ 
and therefore also 
$$
\lim_{n \rightarrow \infty} \Pr[|\Piv_n| > a_n \; | \; \Tr(l_n,k_n) =0 ] =1.
$$
The same argument can be repeated for $-\Tr(l_n,k_n)(-\omega)$. The event that neither $\Tr(l_n,k_n)(\omega)$ nor $\Tr(l_n,k_n)(-\omega)$ happens while the pivotal set of both is larger than $a_n$ still holds with high probability. That is, we find pivotal bits for both $\Tr(l_n,k_n)(\omega)$ and $\Tr(l_n,k_n)(-\omega)$ with high probability and thus push $f_n=\Tr(l_n,k_n)(\omega)- \Tr(l_n,k_n)(-\omega)$ to $1$ or $-1$, respectively.

Furthermore $\Tr(l_n,k_n)(\omega)- \Tr(l_n,k_n)(-\omega)$ is monotone increasing as the sum of monotone increasing functions. 
\epf

Now it only remains to show that with an appropriate choice of the sequences $k_n$ and $l_n$ \eqref{kisVar} and  \eqref{nagyInf} are satisfied.

First, note that
$$
\left(1- \frac{1}{2^{l_n}} \right)^{k_n} \rightarrow 1 \;\text{if and only if} \;  \frac{k_n}{2^{l_n}} \rightarrow 0,
$$
or equivalently
\beq \label{log1Cond} 
\log{k_n} - l_n \rightarrow -\infty,
\eeq
while  after taking the logarithm in both sides \eqref{nagyInf} becomes
\beq \label{log2Cond} 
\log{k_n}+\log{l_n} - l_n \rightarrow \infty.
\eeq
If we now choose $l_n =\log{k_n} + \frac{1}{2}\log{ \log{k_n}}$ then clearly \eqref{log1Cond} is satisfied. As for \eqref{log2Cond}, using that $\log{l_n} \geq \log{ \log{k_n}}$
$$
\log{k_n}+\log{l_n} - l_n \leq \log{k_n} + \log{ \log{k_n}} - (\log{k_n} + \frac{1}{2}\log{ \log{k_n}}) = \frac{1}{2}\log{ \log{k_n}} \rightarrow \infty.
$$
Finally, we note that the argument remains valid with some elementary modifications in case if, instead of the uniform measure we endow the hypercube with the product measure $\Pr_{p}=(1-p\delta_{-1} + p\delta_{1} )^{\otimes k_n }$ for some $p \in (0,1)$. 

\section{Volatility}

Let $X_n(t)$ be the continuous time random walk on the $k_n$ hypercube (where $X_n(0)$ is sampled according to the stationary  measure) with rate $1$ clocks  on the edges. For a sequence of Boolean functions $f_n$  let $C_n$ denote the (random) number of times $f_n(X_n(t))$  changes value in the interval $[0,1]$. The following concepts where introduced in \cite{JS}.

\bde[Volatility, tameness]\label{volat}
A sequence of functions $f_n: \left\{ -1,1 \right\}^{k_n} \longrightarrow  \left\{-1, 1 \right\}$ is called volatile if  the sequence $C_n$ tends to $\infty$ in distribution and  tame, if the sequence $C_n$ is tight.
\ede
It is a (rather intuitive) fact that a non-degenerate noise sensitive sequence is volatile (Proposition 1.17 in \cite{JS}) and all tame sequences are noise stable (Proposition 1.13 in \cite{JS}). The $\Maj$ function is noise stable, but not tame and not volatile either. 

Now we are going to relate our conditions to volatility.
\bl \label{piv_volatile} Let  $f_n: \left\{ -1,1 \right\}^{k_n} \longrightarrow  \left\{-1, 1 \right\}$ be a sequence of Boolean functions with the property that there is a sequence of positive integers $a_n \rightarrow \infty$ such that $\Pr[ |\Piv_n| > a_n] \rightarrow 1$ (where $\Piv_n$ denotes the pivotal set of of $f_n$). Then $f_n$ is volatile. 
\el

\bpf
Let $A_n : = \left\{ |\Piv_n| \leq a_n \right\}$. It is clear that $\E [ \int^{1} _{0}{  \1_{ X_n(t) \in A_n} dt }]  = \Pr[ |\Piv_n| \leq a_n]  \rightarrow 0$ so for every $\epsilon$ for large enough $n$ it holds that
$$
\E [ \int^{1} _{0}{  \1_{ X_n(t) \in A_n} dt }] < \eps^2
$$
and therefore, using Markov's inequality
$$
\Pr [ \int^{1} _{0}{  \1_{ X_n(t) \in A_n} dt } > \eps] < \eps.
$$
By Lemma 1.5 in \cite{JS} volatility is  equivalent with the condition
$$
\lim_n{\Pr[C_n = 0]} = 0.
$$
Now we show that $\Pr[C_n = 0]$ can be arbitrary small. If we choose $n$ large enough so that $e^{- (1-\eps) a_n} < \eps$ 
$$
\Pr[C_n = 0] \leq \Pr [ \int^{1} _{0}{  \1_{ X_n(t) \in A_n} dt } > \eps] +  \Pr [ \int^{1} _{0}{  \1_{ X_n(t) \in A_n} dt } \leq \eps \;\text{and} \; C_n = 0 ] \leq \eps + e^{- (1-\eps) a_n} < 2 \eps,
$$
where we used that $C_n = 0$ can only hold as long as no pivotal bit is switched during the time we are outside of $A_n$.
\epf
We say that the sequences $f_n$ and $g_n$ $o(1)$-close to each other if $\lim_n {\Pr[f_n \neq g_n ]} = 0$. In \cite{Fo} it is proved (Theorem 1.4) that for every sequence of Boolean functions there is a volatile sequence $o(1)$-close to it and in this sense volatile sequences are dense among all sequences of Boolean functions. Our construction has a similar conclusion. Using the fact that  any sequence of Boolean functions can be slightly modified with a bribable sequence in the same way as we did with $\Maj$, we obtain the following strengthening of Theorem 1.4 from \cite{Fo}:

\bc \label{dense_piv}
Any sequence of (monotone) Boolean functions is $o(1)$-close to a (monotone) volatile sequence with the property that $\Pr[\Piv_n > a_n] \rightarrow 1 $ for some sequence of integers $a_n \rightarrow \infty $.
\ec
Although here we consider the uniform measure on the hypercube the same type of questions are meaningful when the uniform measure is replaced by the sequence of product measures $\Pr_{p_n}=(1-p_n\delta_{-1} + p_n\delta_{1} )^{\otimes k_n }$. It has to be noted that Theorem 1.4 in \cite{Fo} is valid for basically all possible sequences $p_n$ under which the question is meaningful, while our construction works in a more restricted range of sequences $p_n$. Most importantly, our results extend to all sequences  $p_n$ that satisfy $0<\liminf{p_n} \leq \limsup{p_n}<1$.  

Furthermore, in \cite{Fo} a  sequence of Boolean functions is constructed which is noise stable and volatile, but at the same time it is not $o(1)$-close to any non-volatile sequence. Such a sequence, of course cannot be obtained with a small modification from some non-volatile stable sequence.

This naturally lead to the following questions:
\begin{quest}
Is there a transitive, noise stable (volatile?) sequence $f_n$ such that  $\Pr[\Piv_n(\omega)\neq \emptyset] \rightarrow 1 $ and $f_n$ is not $o(1)$-close to any sequence which does not have these properties?
\end{quest}
We think that the answer is positive to this question. 

\begin{quest}
Is there a transitive, monotone and noise stable (volatile?) sequence $f_n$ such that  $\Pr[\Piv_n(\omega)\neq \emptyset] \rightarrow 1 $ and $f_n$ is not $o(1)$-close to any sequence which does not have these properties?
\end{quest}
This looks more difficult and it might be the case that the answer is negative.

\section*{Acknowledgements}
Thanks to Christophe Garban, Tom Hutchcroft  and G\'abor Pete for useful comments and discussions. Research supported by ERC Consolidator Grant 772466 “NOISE.

\end{document}